\documentclass[12pt]{article}
\usepackage{graphicx, mathrsfs,amssymb,amsmath,amsfonts,amsbsy,amsthm,latexsym,bm,float}
\begin{document}

\newcommand{\Dc}{\mathcal D}
\newtheorem{teo}{Theorem}[section]
\newtheorem{coro}[teo]{Corollary}
\newtheorem{lem}[teo]{Lemma}

\title{
A simple convergent solver for initial value problems\\
\vskip.5cm}
\author{Rafael G. Campos and Francisco Dom\'{\i}nguez Mota\\
Facultad de Ciencias F\'{\i}sico-Matem\'aticas,\\
Universidad Michoacana, \\
58060, Morelia, Mich., M\'exico.\\
\hbox{\small rcampos@umich.mx, dmota@umich.mx}\\
}
\date{}
\maketitle
{
}\\
\vspace*{1truecm}
\begin{center} Abstract \end{center}
We present a stable and convergent method for solving initial value problems based on the use of differentiation matrices obtained by Lagrange interpolation. This implicit multistep-like method is easy-to-use and performs pretty well in the solution of mildly stiff problems and it can also be applied directly to differential problems in the complex plane.
\vskip1.5cm
\newpage
\section{Introduction}\label{intro}
A general technique for obtaining approximations for the derivative of a well-behaved function $x(t)$ consists in expressing the derivative $x'(t)$ evaluated at an arbitrary point as a linear combination of the function values $x(t_k)$ at $N+1$ nodes $t_k$, $k=0,1,\ldots, N$, i.e., $x'(t)\simeq \sum_{k=0}^N \gamma_k x(t_k)$ where the numbers $\gamma_k$ must be such that the truncation error tends to zero as the mesh goes to zero \cite{Cel92}. This $N+1$-point differentiation scheme can be realized by using the Lagrange interpolating polynomial. In spite of the fact that Lagrange interpolation has been widely and greatly used in numerical differentiation, it seems that a very simple and powerful simplification has been overlooked for many years. If the $N+1$-point derivative of $x(t)$ obtained by Lagrange interpolation is evaluated at the nodes, it can be written the form \cite{Cal83a, Cal83b, Cam04}
\begin{equation}\label{derf}
x'(t_j)=\sum_{k=0}^N \Dc_{jk} x(t_k)+\frac{1}{(N+1)!}P'(t_j)x^{(N+1)}(\tau_j),
\end{equation}
where $\tau_j\in(t_0,t_N)$, $P(t)=\prod_{k=0}^N (t-t_k)$, and
\begin{equation}\label{matdif}
\Dc_{jk}=\begin{cases}\displaystyle \sum_{l \neq j}^N \frac{1}{(t_j-t_l)},&j=k,\\\noalign{\vskip .5truecm}
\displaystyle \frac{P '(t_j)}{(t_j-t_k)\,\,P '(t_k)}, &j\not=k.\\
\end{cases}
\end{equation}
If the function $x(t)$ to be differentiated is a polynomial of degree at most $N$, the truncation error is zero and (\ref{derf}) yield the exact values $x'(t_j)$ from which the function $x'(t)$ can be retrieved by an interpolation. This is why the differentiation matrix $\Dc$ is a projection of $d/dt$ in the subspace of polynomials of degree at most $N$. The arbitrariness of the nodes and the matrix form of Eq. (\ref{derf}), $x'=\Dc x+e$, where $\Dc$ is the matrix whose elements are given by (\ref{matdif}) and $e$ is the vector corresponding to the truncation error, have given rise to a simple method for finding accurate numerical solutions to two-point boundary value problems \cite{Bru90} with a relatively small number of nodes. However, the accuracy achieved by this technique depends on the selection of an auxiliary function needed to incorporate the boundary conditions. On the other hand, the definition (\ref{matdif}) can be generalized straightforwardly to give a differentiation matrix for meromorphic functions, yielding a method for solving singular differential problems in the complex plane \cite{Cam04}.\\
Our aim in this paper is to implement a simple method to solve the initial-value problem by using the differentiation matrix (\ref{matdif}). This is done in the next section. It is shown that this method is convergent and stable and some numerical tests with benchmark problems are given.
\section{Outline of the method}\label{secdos}
Consider the IVP
\begin{equation}\label{ivp}
x'(t)=f(x,t),  \quad x(a)=\alpha,
\end{equation}
defined on ${\mathcal R}=\{(t,x)\vert\,\, a\le t\le b,\,\, c\le x\le d\}$, where $f(x,t)$ is Lipschitz on ${\mathcal R}$ in the variable $x$. Let $a_1$ be a point of  $(a,b)$ and evaluate the differential equation of (\ref{ivp}) at the $N+1$ nodes $a=t_0<t_1<t_2\cdots< t_N=a_1$. This yields the vector equation
\[
\begin{pmatrix} x'(a)\\ x'(t_1)\\ x'(t_2)\\ \vdots\\ x'(t_N)\end{pmatrix}=\begin{pmatrix} f(\alpha,a)\\ f(x_1,t_1)\\ f(x_2,t_2)\\ \vdots\\ f(x_N,t_N)\end{pmatrix}
\]
where $x_k$ stands for $x(t_k)$. According to (\ref{derf}), this equation can be approximated by
\begin{equation}\label{ivpap}
\begin{pmatrix}
\Dc_{aa}& \Dc_{a1}& \Dc_{a2}&\cdots & \Dc_{aN}\\
\Dc_{1a}& \Dc_{11}& \Dc_{12}&\cdots & \Dc_{1N}\\
\Dc_{2a}& \Dc_{21}& \Dc_{22}&\cdots & \Dc_{2N}\\
\vdots&\vdots&\vdots&\ddots&\vdots\\
\Dc_{Na}& \Dc_{N1}& \Dc_{N2}&\cdots & \Dc_{NN}
\end{pmatrix}\begin{pmatrix} \alpha\\ \xi_1\\  \xi_2\\ \vdots\\ \xi_N\end{pmatrix}=
\begin{pmatrix} f(\alpha,a)\\ f(\xi_1,t_1)\\ f(\xi_2,t_2)\\ \vdots\\ f(\xi_N,t_N)\end{pmatrix}.
\end{equation}
This system has only $N$ unknowns $\xi_1, \xi_2,\ldots, \xi_N$, that can be found by solving
\begin{equation}\label{siste}
\sum_{k=1}^N \Dc_{jk} \xi_k-f(\xi_j,t_j)=-\alpha d_j,\quad j=1,2,\ldots,N,
\end{equation}
where $d_j=\Dc_{ja}$. This yields approximations $\xi_j$ to $x(t_j)$ in $[a,a_1]$. Now, let us define $\alpha_1=\xi_N$ and choose $a_2\in (a_1,b]$. Then, the local problem
\[
y'(t)=f(y,t),\quad t\in [a_1,a_2],  \quad y(a_1)=\alpha_1,
\]
can be solved numerically along the above lines to give new approximations $\xi_j$ to $x(t_j)$ in $[a_1,a_2]$. Defining $\alpha_2$ as the new $\xi_N$, this procedure is repeated from the subinterval $[a_{n-1},a_n]$ to the subinterval $[a_n,a_{n+1}]$ until the final point $b$ is reached.\\
\section{Main differences with popular IVP methods}
It must be noted that the proposed method, from now on referred to as the Simple Convergent Solver (SCS), is based on the system of equations (\ref{siste}). Since the solution of the latter produces the approximations to the unknown values $x(t_1), x(t_2),...,x(t_N)$ simultaneously, SCS is an implicit method. It looks like a linear multistep method, however, the fact that $N$ unknowns are calculated simultaneously in every step makes SCS different from those kind of methods. On the other hand, even though SCS is based on Lagrange interpolation and incorporates $N$ values $f(\xi_j,t_j)$, the use of differentiation matrices, instead of quadratures, is what makes SCS also different from the general multistep-multistage methods discussed  in the literature \cite{But03, Sha94}.\\
\section{Analytic properties}
Note that the system (\ref{siste}) has the matrix form 
\begin{equation}\label{ecxin}
D_n\xi_n-f_{\xi_n}=-\alpha_{n-1} d_n,
\end{equation}
where the index $n$ indicates that the variables are being considered in the subinterval $[a_{n-1},a_n]$. In order to simplify the notation, let us rewrite this equation as
\begin{equation}\label{siste20}
D\xi-f_\xi=-\alpha d,
\end{equation}
where $\xi_j=\xi_{nj}$, $(f_\xi)_j=f(\xi_{nj},t_{nj})$, $d_j=d_{nj}=\Dc_{ja}$, and $D=(D_{jk})$ is the $N\times N$ submatrix of $\Dc$ with $D_{jk}\equiv \Dc_{jk}$, $j,k=1,2,\ldots,N$, and $\Dc_{jk}$ being computed according to (\ref{matdif}) at the nodes $a_{n-1}=t_{n0}<t_{n1}<t_{n2}<\ldots< t_{nN}=a_n$.\\
As mentioned above, the approximation to the solution of the IVP (\ref{ivp}) in $[a_{n-1},a_n]$ is obtained by solving the system of equations (\ref{siste20}). In general, for a nonlinear $f(x,t)$, a Newtonian iteration can be applied in order to solve the equations. It must be noted, however, that this is not a major problem: the required Jacobian can be easily approximated, if needed, by using standard differences due to the simple form of the left-hand side, since only the partial derivative $\partial f(x,t)/\partial x$ is required. Thus, the $(k+1)$th iteration can be written as
\begin{equation}
\xi^{(k+1)}=\xi^{(k)}+\eta^{(k)}, \quad (D-\Lambda^{(k)})\eta^{(k)}=-(D\xi^{(k)}-f_\xi^{(k)}+\alpha d)),
\end{equation}
where $\Lambda^{(k)}$ is the diagonal matrix whose non-zero elements are
\[
\begin{pmatrix} {\partial f}(\xi^{(k)}_{n1},t_{n1})/{\partial x}\\ {\partial f}(\xi^{(k)}_{n2},t_{n2})/{\partial x}\\ \vdots\\{\partial f}(\xi^{(k)}_{nN},t_{nN})/{\partial x}\end{pmatrix},
\]
and $f_\xi^{(k)}$ is the vector with elements $(f_\xi^{(k)})_j=f(\xi^{(k)}_{nj},t_{nj})$. The initial trial $\xi^{(0)}$ can be taken as the vector with all its components equal to the initial condition at $a_{n-1}$.\\
On the other hand, if $f$ is linear, let us say $f(x,t)=\kappa x +\phi(t)$, the approximation $\xi$ is the solution of (\ref{siste20}) which becomes 
\[
(D-\kappa 1_N)\xi=-\alpha d-\phi,
\]
where $1_N$ is the identity matrix of dimension $N$ and $\phi$ is the vector of entries $\phi(t_{nj})$, $j=1,2,\ldots,N$.\\
It is clear that the existence of the numerical approximation given by the method is based on the invertibility of the sum of $D$ and a diagonal matrix.\\
The invertibility of $D$ can be easily proved. Let ${\Dc}_N$ be the differentiation matrix constructed according (\ref{matdif}) with the set of $N$ points $t_{n1}<t_{n2}<\ldots< t_{nN}$ (note that the first point $t_{n0}$ has been removed). Then a little algebra shows that
\begin{equation}\label{dtdgar}
(T-t_{n0} 1_N)D=(T-t_{n0} 1_N)\Dc_N+1_N,
\end{equation}
where $T$ is the diagonal matrix whose non-zero entries are $t_{n1},t_{n2},\ldots, t_{nN}$. Since ${\Dc}_N$ is a differentiation matrix, it yields exact values for the derivatives of a polynomial of degree $m$, $m\le N-1$ \cite{Cal83a}. Therefore, the derivatives of $(t-t_{n0})^m$, $m=0,1,\ldots,N-1$, can be reproduced for each $n$ applying ${\Dc}_N$ to the vector of entries $(t_{nj}-t_{n0})^m$, $j=1,2,\ldots,N$. Thus, we have that
\[
(T-t_{n0} 1_N)\Dc_N\begin{pmatrix}(t_{n1}-t_{n0})^m\\ (t_{n2}-t_{n0})^m\\  (t_{n3}-t_{n0})^m\\ \vdots\\ (t_{nN}-t_{n0})^m\end{pmatrix}=
m\begin{pmatrix}(t_{n1}-t_{n0})^m\\ (t_{n2}-t_{n0})^m\\  (t_{n3}-t_{n0})^m\\ \vdots\\ (t_{nN}-t_{n0})^m\end{pmatrix},\quad m=0,1,\cdots,N-1.
\]
The substitution of this result in (\ref{dtdgar}) yields
\[
(T-t_{n0} 1_N)D\begin{pmatrix}(t_{n1}-t_{n0})^m\\ (t_{n2}-t_{n0})^m\\  (t_{n3}-t_{n0})^m\\ \vdots\\ (t_{nN}-t_{n0})^m\end{pmatrix}=
(m+1)\begin{pmatrix}(t_{n1}-t_{n0})^m\\ (t_{n2}-t_{n0})^m\\  (t_{n3}-t_{n0})^m\\ \vdots\\ (t_{nN}-t_{n0})^m\end{pmatrix},\quad  m=0,1,\cdots,N-1.
\]
This equation shows that the eigenvalues of $(T-t_{n0} 1_N)D$ are the integers $1,2,\ldots,N$ for any set of points $t_{n1}<t_{n2}<\ldots< t_{nN}$. Therefore, $(T-t_{n0} 1_N)D$ is invertible. Since $(T-t_{n0} 1_N)$ is invertible, we have proved the
\begin{lem}$D$ is invertible.
\end{lem}
This lemma will be of great importance to prove convergence and stability for the method. More properties of the matrix $D$ will be studied in more detail elsewhere.
\section{Convergence, consistency and stability}
As noted above, the proposed method is different from the standard methods for initial value problems, since it works as a ``block-implicit" method. Thus, it is convenient to discuss {\em ad hoc} proofs of its basic properties.\\
Let us assume that the solution $y(t)$ of the local problem
\[
y'(t)=f(y,t),\quad t\in [a_{n-1},a_n],  \quad y(a_{n-1})=\alpha_{n-1},
\]
is sufficiently smooth and Lipschitz in $y$.\\ 
To prove convergence, let us denote by $y=(y_{n1},y_{n2},\cdots,y_{nN})^T$ the vector whose entries are the values of $y(t)$ at the $N$ nodes $t_{n1},t_{n2},\ldots,t_{nN}$, i.e., $y_{nj}=y(t_{nj})$, and assume that $[a,b]$ has been divided in $M$ subintervals. Then, applying equation (\ref{derf}), we get
\begin{equation}\label{ecerloc}
\sum_{k=1}^N D_{jk} y_{nk}-f(y_{nj},t_{nj})=-\alpha d_j+\frac{E_{nj}}{(N+1)!},\quad j=1,2,\ldots,N,
\end{equation}
where $E_{nj}$ is given by $y^{(N+1)}(\tau_{nj}) P'(t_{nj})$, $\tau_{nj}\in [a_{n-1},a_n]$ and $P(t)=\prod_{k=0}^N (t-t_{nk})$. This can be written in matrix form as
\begin{equation}\label{sistesol}
Dy-f_y=-\alpha d+\frac{E}{(N+1)!},
\end{equation}
where $f_y$ and $E$ are the vectors of entries $(f_y)_j=f(y_{nj},t_{nj})$ and $E_j=E_{nj}$, respectively. From (\ref{siste20}) and (\ref{sistesol}) we get
\[
 \|y-\xi\| \le \|D^{-1}\| \left(\| f_y-f_\xi \|+\frac{\| E\|}{(N+1)!}\right).
\]
Since $f$ is sufficiently smooth and satisfies a Lipschitz condition, there exist constants $\mathcal{L}$ and $\mathcal{K}$ such that
\[
 \|y-\xi\| \le \|D^{-1}\| \left(\mathcal{L}\| y-\xi \|+\frac{\mathcal{K}h^N}{N+1}\right).
\]
where $h=\max_{j=1}^N \vert t_{nj}-t_{n,j-1}\vert$. This yields
\[
 \|y-\xi\| \le\frac{\mathcal{K}\| D^{-1}\|h^N}{(N+1)(1-\mathcal{L}\|D^{-1}\|)}
\]
whenever $\mathcal{L}\| D^{-1}\|<1$. For evenly spaced points with step $h$, $t_{nj} =a_{n-1}+jh$, $j=1,2,\ldots,N$, the matrix $D$ becomes $\Delta/h$, where $\Delta$ is a matrix whose elements are explicitly independent of $h$. Note that $\Delta$ is a matrix whose dimension grows as $h\to 0$. Therefore, the local error between the exact solution $y_n$ and the approximation $\xi_n$ on the $n$th interval becomes
\begin{equation}\label{locerrn}
 \|y_n-\xi_n\| \le\frac{\mathcal{K}_n\| \Delta^{-1}\|h^{N+1}}{(N+1)(1-h\mathcal{L}_n\|\Delta^{-1}\|)}.
\end{equation}
To get a global bound for the true error on $[a,b]$, first let us define $x_n$ as the vector whose entries are the exact values of $x(t)$ at the nodes of $[a_{n-1},a_n]$, i.e., $(x_n)_j=x(t_{nj})$. Now, note that (\ref{locerrn}) gives the true error $\|x_n-\xi_n\|$ whenever the initial condition $\alpha_{n-1}$ (cf. Eq. 
(\ref{ecxin}) for notation) is substituted by the exact value $(x_{n-1})_N$. Let $\delta_n$ be the difference
\[
\delta_n=(x_n)_N-\alpha_n=(x_n)_N-(\xi_n)_N.
\]
Therefore, (\ref{ecerloc}) becomes
\[
\sum_{k=1}^N D_{jk} x_{nk}-f(x_{nj},t_{nj})=-\alpha_{n-1} (d_n)_j-\delta_{n-1} (d_n)_j +\frac{\mathcal{E}_{nj}}{(N+1)!},\quad j=1,2,\ldots,N,
\]
where $\mathcal{E}$ is now the Lagrange interpolation error for $x(t)$. The same argument as above yields
\begin{equation}\label{truerrn}
 \|x_n-\xi_n\| \le \frac{\| \Delta^{-1}\|}{(1-h\mathcal{L}\|\Delta^{-1}\|)}\left(\frac{\mathcal{K}h^{N+1}}{(N+1)}+N h d_M \vert\delta_{n-1}\vert\right).
\end{equation}
where we have used the inequality $\|d_n\|\le N d_M$ and the definition
\[
d_M=\max_{\mathop=\limits_{1\le n\le M}^{1\le j \le N}}\vert (d_n)_j\vert.
\]
By direct inspection we find that for equispaced points
\[
d_M=\frac{P'(t_{n1})}{P'(t_{n0})(t_{n1}-t_{n0})}=\frac{P'(t_{nN})}{P'(t_{n0})(t_{nN}-t_{n0})}=\frac{1}{Nh}.
\]
Taking into account this result and the fact that $\delta_0\equiv 0$, the recursive use of (\ref{truerrn}) yields the global bound for the true error
\[
 \|x_T-\xi_T\| \le \frac{M\mathcal{K}_M\| \Delta^{-1}\|h^{N+1}}{(N+1)(1-h\mathcal{L}_M\|\Delta^{-1}\|)},
\]
where $x_T$ and $\xi_T$ stand for the vectors formed with the exact and approximate values at all the nodes $t_{nj}$, $j=1,\ldots, N$, $n=1,\dots, M$.  This proves the
\begin{teo}
The method given by (\ref{siste20}) is convergent and of order $\mathcal{O}(h^{N})$.
\end{teo}
\noindent 
In addition, it is important to note that, as discussed by Celia and Gray \cite{Cel92}, a difference approximation as that given by the differentiation matrix (\ref{ivpap}) and from which SCS is derived, satisfies the equation
\[
\sum_{j=0}^N \gamma_{kj}\, x(t_j)=\frac{dx(t_k)}{dt}+\mathcal{O}(h^p),\quad 0\le k\le N,\quad p>0, 
\]
imposed by the consistency requirement, if
\[
\begin{pmatrix}
1& 1& \cdots & 1\\
(t_0-t_k)& (t_1-t_k)& \cdots & (t_N-t_k)\\
\end{pmatrix}\begin{pmatrix} \gamma_{k1}\\ \gamma_{k2}\\  \gamma_{k4}\\ \vdots\\ \gamma_{kN}\end{pmatrix}=
\begin{pmatrix} 0\\  1\end{pmatrix}
\]
which is a straightforward consequence of the properties of the differentiation matrix (\ref{matdif}). The analogous expression for higher order derivatives can be obtained in a similar manner. Therefore, the proposed method, SCS, is consistent.\\
Furthermore, conditional stability of the difference scheme with respect to small changes on the initial value can also be easily proved. Let $\alpha$ and $\tilde{\alpha}$ be different initial values for a subinterval, and $\xi_\alpha$, $\xi_{\tilde{\alpha}}$ the corresponding solutions calculated with the method. Then
\[
\|\xi_\alpha-\xi_{\tilde{\alpha}}\|\le\frac{\mathcal{K}\|D^{-1}\|\,|\alpha-\tilde{\alpha}|}{1-\mathcal{L}\|D^{-1}\|}.
\]
Therefore, small changes in the initial condition yield small changes in the numerical solution whenever $\mathcal{L}\| D^{-1}\|<1$.
\section{Numerical tests}
In this section we test SCS with some benchmark problems used frequently in the literature. SCS is tested against two quality MATLAB solvers, ODE15s, with default parameters (order 5) and quasi-constant step size and ODE45 \cite{Sha97, Sha03}. It should be noted that SCS, as used here, is not based on error control algorithms or adaptive step size strategies. We consider default blocks with $N=5$ equispaced points in each subinterval (in order to compare methods of the same order) until the whole time interval is covered. In the following tables SCS stands for the proposed method and $E$ is the absolute error and $\| E \|$ stands for the euclidean norm for one-dimensional problems or the Frobenius norm for two-dimensional problems.\\
All the tests were run in MATLAB 7.2 using a personal computer with 2Gb RAM and Intel $\copyright$ processor running at 1.99GHz.
\subsection{First Order equations}
\subsubsection*{Example 1}
Let us consider the stiff problem \cite{Sha03}
\begin {equation}\label{pviuno}
x'=-100x+10,\quad x(0)=1,
\end{equation}
for $0\le x \le 0.2$. The solution is $x(t)=(1+9\,e^{-100\,t})/10$ and the results are shown in Table 1 . They are compared with those obtained with ODE15s at some points of $[0,0.2]$
\begin{table} 
\begin{center}
\begin{minipage}{12cm}
\caption
{\small Results for the IVP (\ref{pviuno}) at selected points of $[0,0.2]$. $E_{\text{SCS}}$ and $E_{\text{ODE15s}}$ are the absolute errors obtained by using SCS and ODE15s respectively. The corresponding norms are $\| E_{\text{SCS}} \|=  0.0000714$ and $\| E_{\text{ODE15s}} \|= 0.000242$.}
\end{minipage}
\vskip.2cm
\label{tab1}\begin{tabular}{rrr}
\hline\noalign{\smallskip}
\noalign{\smallskip}
$\hfill t\hfill$&$\hfill E_{\text{SCS}}\hfill$ &$\hfill E_{\text{ODE15s}}\hfill $\\
\noalign{\smallskip}\hline\noalign{\smallskip}
0.00 & 0.0                        & 0.0\\
0.02 & 0.0000688546 & 0.000219504\\
0.04 & 0.0000186422 & 0.0000679712\\
0.06 & 3.78549$\times10^{-6}$ & 0.0000572691\\
0.08 & 6.83273$\times10^{-7}$& 0.0000134028\\
0.10 & 1.15621$\times10^{-7}$& 0.0000405086\\
0.12 & 1.87825$\times10^{-8}$& 4.75526$\times^{-6}$\\
0.14 & 2.96643$\times10^{-9}$& 2.15198$\times^{-7}$\\
0.16 & 4.5894$\times10^{-10}$ & 0.0000226791\\
0.18 & 6.9895$\times10^{-11}$& 0.0000164057\\
0.20 & 1.0513$\times10^{-11}$& 3.69695$\times10^{-6}$\\
\noalign{\smallskip}\hline
\end{tabular}
\end{center}
\end{table}
\subsubsection*{Example 2}
Let us now consider the simple problem 
\begin {equation}\label{pvidos}
x'=100x,\quad x(0)=1,
\end{equation}
for $0\le x \le 0.1$. Here, $x(t)=e^{100 t}$, and the results are shown in Table 2 and compared with those obtained with ODE15s at some points of $[0,0.1]$
\begin{table}\label{tab2} 
\begin{center}
\begin{minipage}{12cm}
\caption
{\small Results for the IVP (\ref{pvidos}) at selected points of $[0,0.1]$. $E_{\text{SCS}}$ and $E_{\text{ODE15s}}$ are the absolute errors obtained by using SCS and ODE15s respectively. The corresponding norms are $\| E_{\text{SCS}} \|=8.03$ and $\| E_{\text{ODE15s}} \|=396.2$.}
\end{minipage}
\vskip.2cm
\begin{tabular}{rrr}
\hline\noalign{\smallskip}
\noalign{\smallskip}
$\hfill t\hfill$&\hfill $E_{\text{SCS}}$($\times 10^2$)\hfill &\hfill $E_{\text{ODE15s}}$($\times 10^2$)\hfill \\
\noalign{\smallskip}\hline\noalign{\smallskip}
0.00 & 0.0                        & 0.0\\
0.02 & 0.00000535 & 0.00046360\\
0.04 & 0.00007917 & 0.00452663\\
0.06 & 0.00087755 & 0.04615606\\
0.08 & 0.00864604& 0.44160702\\
0.10 & 0.07986052& 3.93706426\\
\noalign{\smallskip}\hline
\end{tabular}
\end{center}
\end{table}
\subsubsection*{Example 3}
Next, we consider the nonlinear problem \cite{Bur00}
\begin {equation}\label{pvitres}
x'=5e^{5t}(x-t)^2+1,\quad x(1)=-1,
\end{equation}
for $0\le t \le 1$. For this equation, $x=y-e^{-5t}$. The results are shown in Table 3 and compared with those obtained with ODE45 and ODE15s at some points of $[0,1]$. \\
\begin{table}\label{tab3}
\begin{center}
\begin{minipage}{12cm}
\caption
{\small Results for the IVP (\ref{pvitres}) at selected points of $[0,1]$. $E_{\text{SCS}}$, $E_{\text{ODE15s}}$ and $E_{\text{ODE45}}$ are the absolute errors obtained by using SCS, ODE15s and ODE45 respectively. The corresponding norms are $\| E_{\text{SCS}} \|=6.7\times 10^{-9}$, $\| E_{\text{ODE15s}} \|=6\times 10^{-4}$ and $\| E_{\text{ODE45}} \|=3\times 10^{-4}$.}
\end{minipage}
\vskip.2cm
\begin{tabular}{rrrr}
\hline\noalign{\smallskip}
\noalign{\smallskip}
$\hfill t\hfill$&$\hfill E_{\text{SCS}}\hfill$ &$\hfill E_{\text{ODE15s}} \hfill$&$\hfill E_{\text{ODE45}} \hfill$\\
\noalign{\smallskip}\hline\noalign{\smallskip}
0.2	&	5.19952E-10	&	5.74683E-04	&	2.32239E-04	\\
0.4	&	6.99985E-11	&	7.58958E-05	&	1.74210E-04	\\
0.6	&	9.39138E-12	&	1.69276E-04	&	7.74205E-05	\\
0.8	&	1.13487E-12	&	6.44637E-05	&	3.22118E-05	\\
1.0	&	6.68797E-09	&	4.02615E-05	&	1.07195E-05	\\
\noalign{\smallskip}\hline
\end{tabular}
\end{center}
\end{table}
\subsection{Second order problems}
In this subsection, we compare the performance of SCS against ODE15s and ODE45 in second order problems.
\subsubsection*{Example 4}
Let us consider the problem \cite{Sha94}
\begin {equation}\label{pvicuatro}
x'_1=-0.1 x_1-199.9 x_2, \qquad x'_2=-200 x_2,
\end{equation}
with $x_1(0)=2$, $x_2(0)=1$, $0\le t \le 50$. The solution is given by
\[
x_1(t)=\exp(-0.1t)+\exp(-200t)\qquad x_2(t)= \exp(-200t),
\]
and the results corresponding to the first component of the solution at some points are shown in Table 4.
\begin{table}\label{tab4}
\begin{center}
\begin{minipage}{12cm}
\caption
{\small Results for the IVP (\ref{pvicuatro}) at selected points of $[0,50]$. $E_{\text{SCS}}$, $E_{\text{ODE15s}}$ and $E_{\text{ODE45}}$ are the absolute errors obtained by using SCS, ODE15s and ODE45 respectively. The corresponding Frobenius norms are $\| E_{\text{SCS}} \|=1.1256\times 10^{-3}$, $\| E_{\text{ODE15s}} \|=1.7566$ and $\| E_{\text{ODE45}} \|=1.7558$.}
\end{minipage}
\vskip.2cm
\begin{tabular}{rrrr}
\hline\noalign{\smallskip}
\noalign{\smallskip}
$\hfill t\hfill$&$\hfill E_{\text{SCS}}\hfill$ &$\hfill E_{\text{ODE15s}} \hfill$&$\hfill E_{\text{ODE45}} \hfill$\\
\noalign{\smallskip}\hline\noalign{\smallskip}
10.00	&	4.35870E-04	&	3.67880E-01	&	3.67880E-01	\\
20.00	&	4.32250E-05	&	1.35340E-01	&	1.35340E-01	\\
30.00	&	2.37190E-05	&	4.97870E-02	&	4.97870E-02	\\
40.00	&	1.16350E-05	&	1.83160E-02	&	1.83150E-02	\\
50.00	&	5.35100E-06	&	6.73790E-03	&	6.73720E-03	\\
\noalign{\smallskip}\hline
\end{tabular}
\end{center}
\end{table}
\subsubsection*{Example 5}
Let us consider the Lotka-Volterra system
\begin{equation}\label{pvicinco}
x'_1=x_1\,(0.76-0.45\,x_2), \qquad x'_2=-x_2\,(0.18-0.82\,x_1),
\end{equation}
with $x_1(0)=0.1$, $x_2(0)=0.1$, $0\le t \le 1$. The solver ODE45 requires 67 function evaluations in 11 successful steps. SCS requires 36 block function evaluations and requires to solve 5 linear systems. The results can be compared in Table 5.\\
\begin{table}\label{tab5}
\begin{center}
\begin{minipage}{12cm}
\caption
{\small Differences $\vert (x_1)_{\text{SCS}}- (x_1)_{\text{ODE45}}\vert$ and $\vert (x_2)_{\text{SCS}}-(x_2)_{\text{ODE45}}\vert$
of the components $(x_1(t),x_2(t))$ for the LotkaVolterra system (\ref{pvicinco}), computed by using SCS and ODE45 respectively at selected points of $[0,1]$.}
\end{minipage}
\vskip.2cm
\begin{tabular}{ccc}
\hline\noalign{\smallskip}
\noalign{\smallskip}
$\hfill t\hfill$&$\hfill \vert (x_1)_{\text{SCS}}- (x_1)_{\text{ODE45}}\vert\hfill$ &$\hfill \vert (x_2)_{\text{SCS}}-(x_2)_{\text{ODE45}}\vert \hfill$\\
\noalign{\smallskip}\hline\noalign{\smallskip}
0.25	&	7.13490E-09	&	1.18070E-09	\\
0.50	&	1.68620E-08	&	2.97240E-09	\\
0.75	&	2.93810E-08	&	5.61470E-09	\\
1.00	&	4.54880E-08	&	9.59720E-09	\\
\noalign{\smallskip}\hline
\end{tabular}
\end{center}
\end{table}
\section{Final Remarks}\label{secseis}
The cornerstones of a simple and convergent method to solve initial value problems in ordinary differential equations have been presented. Due to its simplicity, it can be easily implemented without the need of knowing previous function values except for the initial value. The implementation can be made with equispaced points (as in this work) or can be made adaptive.\\
The numerical tests show that it is indeed a competitive option, which produces accurate results in a wide range of problems. Only some very specific examples were discussed in this paper, but SCS can be easily extended to solve high order, non linear, and vector systems of differential equations. It is important to note that SCS is not only a competitive method, but also a innovative one, since it is a {\em block implicit} method, which makes it different form the standard implicit methods for the numerical solution of ordinary differential equations.
In future papers, the application of these ideas to boundary problems and partial differential equations will be discussed.


\end{document}